\documentclass[12pt]{article}
\usepackage{amsmath}
\usepackage{amsfonts}
\usepackage{amssymb}
\usepackage{amsthm}
\usepackage[all]{xy}
\usepackage[usenames,dvipsnames]{color}
\usepackage{hyperref}
\newcommand{\C}{\mathbb C}
\newcommand{\N}{\mathbb N}
\newcommand{\F}{\mathbb F}
\newcommand{\wP}{w\mathbb P}
\newcommand{\A}{\mathbb A}
\newcommand{\Q}{\mathbb Q}
\newcommand{\Z}{\mathbb Z}
\newcommand{\cH}{\mathcal H}
\newcommand{\broken}{\dasharrow}
\DeclareMathOperator{\Grass}{Grass}
\DeclareMathOperator{\Pic}{Pic}

\DeclareMathOperator{\Spec}{Spec}
\DeclareMathOperator{\Bir}{Bir}
\DeclareMathOperator{\Tr}{Tr}
\newcommand{\magma}{{\sc Magma}}
\renewcommand\theta{\vartheta}
\renewcommand\P{{\mathbb P}}
\renewcommand\phi{{\varphi}}
\newcommand\eps{{\varepsilon}}
\newcommand\kbar{{\overline{k}}}
\newcommand\Xbar{{\overline{X}}}
\DeclareMathOperator{\mult}{mult}
\newtheorem{thm}{Theorem}

\newtheorem{lem}[thm]{Lemma}

\theoremstyle{definition}

\theoremstyle{remark}

\newtheorem*{rem*}{Remark}
\theoremstyle{remark}

\numberwithin{equation}{section}
\numberwithin{thm}{section}
\begin{document}
\title{Computational birational geometry\\of minimal rational surfaces}
\author{Gavin Brown \and Alexander Kasprzyk \and Daniel Ryder}
\date{}
\maketitle
\begin{abstract}
  The classification of minimal rational surfaces and the
  birational links between them by Iskovskikh, Manin and others
  is a well-known subject in the theory of algebraic surfaces.
  We explain algorithms that realise links of type~II between
  minimal del Pezzo surfaces, one of the major classes of
  birational links, and we describe briefly how this fits into a
  large project to implement the results of Iskovskikh's
  programme in \magma.
\end{abstract}

The theory of minimal rational surfaces and their birational links
works over an arbitrary perfect field~$k$.
Our interest here is the case $k=\Q$ or a number field, in part because
we can compute in these fields, but also because the implementation
we present here imposes some conditions on the characteristic.
We let $\kbar$ denote
the algebraic closure of~$k$. When a geometric object $A$ is
defined over~$k$, $\overline{A}$ denotes its base change to
$\kbar$.  All algebraic surfaces in this paper are nonsingular.
A surface $X$ is {\em rational}\/ if there is an isomorphism
$\Xbar\buildrel{\cong}\over{\longrightarrow}\P^2_\kbar$ defined
over~$\kbar$ (but not necessarily over~$k$).  This notion is
sometimes also called {\em geometrically rational}, to emphasise the
algebraic closure in the statement.  A surface $X$ is {\em minimal}\/ if
any birational morphism $X\rightarrow Y$ defined
over~$k$ is an isomorphism.

When $k=\C$ the results of the theory are over a century old.
The use of Mori theory and the Sarkisov programme, following
Corti \cite{C} and Iskovskikh \cite{I}, make our approach here
very different to what it would have been 15 years ago. Yet more
recent results, such as Hacon and McKernan's approach to the
Sarkisov programme, may change it again in the future.

As an introduction, we outline the results of the theory of
minimal rational surfaces.  This theory is the end result of 150
years of development: from Cayley's computation \cite{Ca} of the
27 lines on a cubic surface around 1850 (in effect, computing the
divisor class group); via Castelnuovo's rationality criterion of
around 1900 (which determines, over an algebraically closed
field, whether or not the surface can be rationally
parametrised); Segre's analysis \cite{Se} in 1942 of the
nonrationality of a cubic surface in terms of its Picard group,
working over an arbitrary field; and Manin \cite{M} and
Swinnerton-Dyer's \cite{SD} geometrical analysis of rationality
in the 1960s in terms of points and configurations of the 27
lines; and culminating in Iskovskikh's complete analysis in a
series of papers during the 1980s.  A full account of the modern
theory is given in~\cite{I}; it is also sketched in the appendix
to~\cite{C}, together with details of the \hbox{(uni)rationality}
of these surfaces over the ground field.  The book
\cite{KSC}~gives proofs of several parts of the theory.

We sketch the classical and modern theories in Section~\ref{sec:1}.
Our main technical tools are explained in Section~\ref{sec:2}, and
they are applied in Section~\ref{sec:3} to compute the major
birational maps of the theory: Geiser and Bertini involutions.
In Section~\ref{sec:group} we show how our algorithms can be
applied to analyse the Cremona group of birational selfmaps of a
minimal cubic surface. Section~\ref{sec:future} explains how
the work described here fits into a broader research programme.

We have implemented the algorithms described below in the
computational algebra system \magma~\cite{BCP};
all code and examples are available at \cite{BKR}.

\section{Minimal models of rational surfaces}\label{sec:1}
\subsection{Classical results over the complex numbers}\label{sec:classical}
A complete analysis of this case is in \cite{KSC} Chapters~2
and~3.  Over~$\C$, a surface is minimal if and only if it does
not contain any curves $C\cong\P^1$ with $C^2=-1$.

The minimal rational surfaces over~$\C$ are: $\P^2$, the
minimal rational surface scrolls $\F_n$ for $n\ge 2$, and
$\F_0=\P^1\times\P^1$, which is often embedded as a quadric
$(xy=zt)\subset\P^3$.  \magma\ treats scrolls as ambient spaces
using multigraded rings:\footnote{This is not in the 2008 export
  of \magma, but will be included in later versions.  More
  generally, any implementation of toric geometry contains surface
  scrolls among its first examples.}  the scroll $\F_n$ has
homogeneous coordinate ring $k[u,v,x,y]$ bi-graded by the columns
of the matrix
\[
\left(
\begin{array}{cc|cc}
1&1&0&-n\\0&0&1&1
\end{array}
\right)
\]
and irrelevant ideal $(u,v)\cap(x,y)$, which is indicated by the
separating vertical line in the matrix.  From the point of view
of Mori theory, it is crucial to consider each surface with a
chosen structure map: for $\P^2$ it is the trivial map to a
point, and for the scrolls $\F_n$ with $n\ge 2$ it is their
natural map to $\P^1$.  The quadric $\F_0$ admits two maps to
$\P^1$, the two projections to its factors, and we must choose
one of these.
Although we often omit the map in our notation,
it is illiterate not to have it (or, equivalently, the
corresponding choice of extremal ray) in mind.

We describe four classes of \emph{birational links\/}, or
elementary birational maps, between these surfaces (including
$\F_1$, although it is not minimal):
\begin{description}
\item{I} (blowup) $\P^2\broken\F_1$ in suitable coordinates by
  $(u,v,w)\mapsto (u,v,w,1)$.
\item{II} (elementary transformation) either $\F_i\broken
  \F_{i-1}$ in coordinates by $(u,v,x,y)\mapsto (u,v,x,uy)$, or
  $\F_i\broken \F_{i+1}$ in coordinates by $(u,v,x,y)\mapsto
  (u,v,ux,y)$.
\item{III} (blowdown) $\F_1\rightarrow\P^2$ in suitable
  coordinates by $(u,v,x,y)\mapsto (uy,vy,x)$.
\item{IV} (change factors) the identity map on $\F_0$, but taking
  the two different projection maps $\F_0\rightarrow \P^1$ on the
  source and target.
\end{description}
The main result is that any birational map between minimal
rational surfaces factorises as a composition of these elementary
birational maps.
\begin{thm}[Noether--Castelnuovo]\label{thm:classical}
  Let $X$ and~$Y$ be minimal rational surfaces over~$k=\C$ and
  $\phi\colon X\broken Y$ a birational map between them.  Then
  there are birational links $\eps_1,\dots,\eps_r$ and an
  automorphism $\psi$ of $Y$ such that
  $\phi = \psi\circ\eps_r\circ\cdots\circ\eps_1$.
\end{thm}

For example, if $X=Y=\P^2$ and $\phi\colon(x,y,z) \mapsto
(1/x,1/y,1/z)$ is the standard quadratic Cremona transformation
then $\phi$ factorises as
\[
\P^2 \broken \F_1 \broken \F_0 \broken \F_1 \rightarrow \P^2.
\]
Indeed, as a composition of maps, up to a linear automorphism of
$\P^2$ it is
\[
\begin{array}{lll}
\multicolumn{2}{l}{(x,y,z) \mapsto (x,y,z,1) \mapsto (x,y,z,x)
  \mapsto} &\\
\qquad & \multicolumn{2}{r}{(x,y,z,xy) = (xz,yz,1,xy) \mapsto
  (xz,yz,xy) = (1/y,1/x,1/z),}
\end{array}
\]
corresponding to the sequence: blowup in $(0,0,1)$, elementary
transform in $(0,1,0,1)$, elementary transform in $(1,0,0,1)$,
blowdown of the negative section. 
\subsection{Minimal rational surfaces over a perfect field}\label{sec:minratsfs}
We will be concerned with two classes of surfaces.
First, the class of del Pezzo
surfaces of Picard rank~1; such surfaces are automatically
minimal since any nontrivial (birational) morphism $X\rightarrow
Y$ decreases the Picard rank.  Second, the class of conic bundles
over a smooth rational curve. Surfaces in this class are not
necessarily minimal; the surface $\F_1$ is not minimal, for instance.
Nevertheless, a minimal rational surface over~$k$
belongs to one of these two classes.

\paragraph{Minimal del Pezzo surfaces.}
These are surfaces $X$ with $-K_X$ ample and $\Pic(X)\cong\Z$.  In
their \hbox{(pluri-)anticanonical} embedding, they are in one of
the families listed in Table~\ref{10dP};
in this embedding, the divisor of any degree~1 linear
section of $X$ is linearly equivalent to $-K_X$. 
\begin{table}[ht]
\centering
$
\begin{array}{cl}
  \text{\ $K_X^2$\ } & \text{Description} \\
  \hline
  1 & X_6 \subset\P(1,1,2,3) \\
  2 & X_4 \subset\P(1,1,1,2) \\
  3 & X_3 \subset\P^3	\\
  4 & X_{2,2} \subset\P^4	\\
  5 & X \subset\P^5 \text{ a general linear section of
    $\Grass(2,5)\subset\P^9$}	\\
  6 & X\subset\P^6 \text{ defined by 9 quadrics} \\
8 & X\subset\P^8   \text{ defined by 20 quadrics}  \\
9 & X\subset\P^9 \text{ defined by 27 quadrics}
\end{array}
$
\caption{Families containing the minimal del Pezzo surfaces. The case $K^2=8$
splits into two, according to whether $X$ is isomorphic to a quadric
$X\cong X_2\subset\P^3$ or not.
The case $K^2=9$
splits into two, according to whether $X\cong\P^2$ or not.}
\label{10dP}
\end{table}
In the table, the entries $X\subset \P^5$, $X\subset \P^6$, $X\subset
\P^8$, $X\subset \P^9$ are surfaces of degree~$d$ in~$\P^d$; these
are not our main concern here, so we give only coarse information
about their equations (although see Section~\ref{sec:rat} for more results).
In degree $d=5$, the ideal of
$X\subset\P^5$ is generated by the five maximal pfaffians of a
skew $5\times 5$ matrix of linear forms---that is one
interpretation of the role of the Grassmannian.

We emphasise
that whether a given surface~$X$ in Table~\ref{10dP} actually
belongs to the class we are discussing depends on whether
$\Pic(X) \cong \Z$; this, in turn, depends on both~$k$ and on the
defining equations of~$X$, and in general it is a difficult problem.
If so then
the Picard group is generated by the class of $-K_X$ in all cases
except the
nonsingular quadric $X_2\subset\P^3$ (a special case for $K^2=8$),
when $-\frac{1}{2}K_X$ generates,
and
$\P^2$ (a special case for $K^2=9$), when
$-\frac{1}{3}K_X$ generates.

\paragraph{Conic bundles.}
These are surfaces $X$ that admit a morphism $f\colon X
\rightarrow C$ to a smooth rational curve $C$ such that $\Pic(X)
= f^*\Pic(C) + \Z$ and a general geometric fibre $F$ of~$f$ is a
conic $F_2\subset\P^2$. The base curve $C$ may have $k$-rational
points (in which case it is isomorphic to $\P^1$) or not (in
which case it is itself isomorphic to a plane conic).

If $q_1,\dots,q_m\in C$ are the closed points at which $f$ is
degenerate---the irregular values of~$f$---then the degree of~$X$
is
\[
K_X^2 = 8 - \deg(q_1) - \cdots - \deg(q_m).
\]

When $C\cong\P^1$, the surface $X$ can be written as a relative
anticanonical model:
\[
X  = (F=0)
\quad\text{in the scroll}\quad
\F_{a,b} =
\left(
\begin{array}{cc|ccc}
1&1&0&-a&-b\\
0&0&1&1&1
\end{array}
\right)
\]
where $0\le a\le b$ and $d\ge 0$ are all integers, and $F=F(u,v,x,y,z)$,
in coordinates $u,v,x,y,z$ on $\F_{a,b}$ of bi-degrees given by the columns
of the matrix,
is bi-homogeneous of weight $(d,2)$;
in particular, $F$ is a quadric in $x,y,z$.
(Although only $d\ge -a$ is required for this linear system to contain
irreducible surfaces, such $X$ has a section whenever $d < 0$ and
so is not minimal unless $X\cong\F_n$;
compare with Lecture~2 of \cite{R}.)
In this notation, it is easy to compute
\[
K_X^2 = 8 - 3d - 2(a+b),
\]
so $\sum\deg(q_i) = 3d + 2(a+b)$.
Suitable diagonal surfaces
\[
X = (A_dx^2+B_{d+2a}y^2 + C_{d+2b}z^2=0) \subset \F_{a,b}
\]
with $A, B, C$ forms in $u,v$ of the indicated degrees give
examples with arbitrarily complicated irregular values
(the roots of $ABC=0$ where the remaining quadric is irreducible over~$k$).

\subsection{Birational links}\label{sec:links}
The surfaces we consider are examples of two-dimensional Mori fibre spaces,
that is, maps $f\colon X\rightarrow S$ with $X$ a surface, $S$ a point
or a nonsingular curve and $f$ a morphism with connected fibres,
$-K_X$ relatively ample and of relative Picard rank~$1$. The map $f$ is simply
the given map for a conic bundle, and it is the trivial map to the point $\Spec k$,
denoted also by $\{ * \}$ with $k$ implicit,
when $X$ is a minimal del Pezzo surface.

A \emph{birational link\/} between \hbox{two-dimensional} Mori
fibre spaces $f\colon X\rightarrow S$ and $f'\colon X'\rightarrow
S'$ is a diagram
\begin{equation} \label{eqn:link}
\xymatrix@C=0.8cm{
  X \ar[d]_f \ar@{-->}[rr]^F & & X \ar[d]^{f'} \\
  S & & S' }
\end{equation}
in which $F$ is a birational map arising in one of four ways:
\begin{description}
\item[Type I.]
These are commutative diagrams of the form
\[
\xymatrix@C=0.8cm{
  X \ar[d] & & X' \ar[ll]_{F^{-1}} \ar[d] \\
  S & & S' \ar[ll]} \] where $F^{-1}$ is the blowup of an
irreducible closed point of~$X$.  For example, take $X = \P^2$,
$S = \{*\}$ (a point), $X' = \F_1$ and $S' = \P^1$.

\item[Type II.]  In this case there is a surface $Y$ and maps
  $h\colon Y\rightarrow X$ and $h'\colon Y\rightarrow X'$ fitting
  into a commutative diagram:
\[
\xymatrix@C=0.8cm{
  & Y \ar[dl]_h \ar[dr]^{h'} \\
  X \ar[d] & & X' \ar[d] \\
  S \ar[rr]^{\cong}& & S' } \]
The maps $h$ and~$h'$ are the blowups of irreducible closed
points of $X$ and~$X'$ respectively.
\item[Type III.]  These are inverses of links of type~I, so there
  is a commutative diagram
\[
\xymatrix@C=0.8cm{
  X \ar[rr]^{F} \ar[d] & & X' \ar[d] \\
  S \ar[rr] & & S' } \] in which $F$ is the blowup of an
irreducible closed point of~$X'$.

\item[Type IV.] Here we have a diagram of the form
\[
\xymatrix@C=0.8cm{
  X \ar@{=}[rr] \ar[d] & & X' \ar[d] \\
  S \ar[dr] & & S' \ar[dl] \\
  & \Spec k & } \] in which $X$ and~$X'$ are the same surface but
the link changes the Mori fibre space structure.  For example,
take $X = \F_0 = \P^1 \times \P^1$ with $X \to S$ the projection
onto the first factor~$\P^1$ and $X = X' \to S'$ the projection
onto the second factor.

\end{description}

\medskip\noindent
The following central result is the analogue of
Theorem~\ref{thm:classical}.

\begin{thm}[Iskovskikh \cite{I} Theorem 2.5] \label{thm:decomp}
  Let $X$ and~$Y$ be minimal rational surfaces over a perfect
  field~$k$ and $\phi\colon X\broken Y$ a birational map between
  them.  Then there are birational links $\eps_1,\dots,\eps_r$
  and an automorphism $\psi$ of~$Y$, all defined over~$k$,
  such that $\phi = \psi\circ\eps_r\circ\cdots\circ\eps_1$.
\end{thm}

\subsection{Geiser and Bertini involutions}
In addition to the factorisation theorem above, Iskovskikh \cite{I} Theorem~2.6
classifies all the birational links~\eqref{eqn:link} that can occur
into 41 different classes according to the type of $X$ and~$X'$.
In Sections~\ref{sec:2}--\ref{sec:3}, we describe the tools central
to the implementation of seven of these classes:
the Geiser and Bertini involutions on
del Pezzo surfaces of degrees~2, 3, 4 and~5.
These are all links of type~II.

We use the following notation throughout. The degree~1 hyperplane
section of $X\subset\P^d$ is denoted by $A$, and $|nA|$ denotes
the linear system of all sections of degree~$n$. 
We regard elements of $|nA|$ both as homogeneous polynomials on
the whole of $\P^d$ and as the zero loci on $X$ they define, and
we move freely between these two descriptions. 
(When $n$ is large, we
may take a vector space complement to those polynomials that lie
in the ideal of~$X$---they do not define divisors on $X$, of course.)
A subsystem is denoted by $\cH\subset |nA|$; this is
some linear subspace of polynomials of degree~$n$.
If $P\in X$ is a point, then the space of polynomials of degree~$n$
whose zero
loci on $X$ vanish to order $m$ at $P$ is denoted by $|nA-mP|$
(or $\cH(-mP)$ if restricting attention to a subsystem $\cH$).

\subsubsection*{Geiser involutions}
Let $X\subset\P^d$ be a del Pezzo surface of degree $d=3$, $4$
or~$5$ and $P\in X$ a closed point of degree~$d-2$.  The
\emph{Geiser involution $i_P\colon X \broken X$ with centre~$P$}
is defined as follows:  $P$~spans a linear $\Pi=\P^{d-3}$, and a
general hyperplane $H=\P^{d-2}$ containing $\Pi$ intersects $X$
in $d$~points, the sum of~$P$ and an effective 0-cycle~$Q$ of
degree~2; $i_P$~exchanges the geometric points of~$Q$ (whether
they are defined over~$k$ or not).

This clearly defines an involution of~$X$ and it is
straightforward to check that $i_P$ is a rational, and hence
birational, map.  Moreover, if we suppose that $X$ is minimal,
following Iskovskikh the linear system corresponding to~$i_P$
is~$|(d-1)A - dP|$.

\begin{rem*}
  The minimality condition may be surprising here; the point is
  that if $(E \subset Y) \to (P \in X)$ is the blowup of $P \in
  X$ then Riemann--Roch computes \[ \chi \big( (d-1)A - dE \big) =
  d+1, \]but we need to prove that this is
  equal to $h^0\big( (d-1)A - dE \big)$, which requires the
  first cohomology to vanish.  Minimality ensures this, although
  in practice it is often the case for
  nonminimal examples too.  The same remark holds for Bertini
  involutions.
\end{rem*}

Now we describe $i_P$ as a link of type~II, in the style
of Section~\ref{sec:links}.  We form the diagram
\[
\xymatrix@C=0.8cm{
  & Y \ar[dl]_f \ar[dr]^g \\
  X \ar[d] & & X \ar[d] \\
  \{*\} & & \{*\} } \] in which $f \colon Y \to X$ is the blowup
of~$P$, the Mori fibre space structure $X \to \{*\}$ on~$X$ is
the trivial one, mapping to a point, and $g \colon Y \to X$ is a
morphism that we define below.

It can be shown that there is a unique effective curve $C$ on~$X$
that is of anticanonical degree $d-2$ and has multiplicity $d-1$
at each geometric point $P' \in P$: for example, if $d=3$ then
$C$ is the tangent curve $T_P(X) \cap X$.  Let $F \subset Y$ be
the strict transform of~$C$.  Then the following hold (we omit
the proofs): $F$ is a $(-1)$-curve or a union of conjugate
$(-1)$-curves, its contraction~$g$ maps to~$\P^d$, and for some
choice of~$g$---which is only defined up to an automorphism
of~$\P^d$ by what we have said so far---we have $g \colon Y \to X
\subset \P^d$ and the composite $g \circ f^{-1}$ is equal to~$i_P
\colon X \broken X$.

\subsubsection*{Bertini involutions}
Let $X\subset\P^d$ be a del Pezzo surface of degree $d=2$, 3, 4
or~$5$ and $P\in X$ a closed point of degree $d-1$.  The
\emph{Bertini involution $j_P\colon X \broken X$ with centre~$P$}
is defined as follows:  $P$~spans a linear $\Pi=\P^{d-2}$, and a
general hyperplane $H=\P^{d-1}$ containing $\Pi$ intersects $X$
in a nonsingular curve $C$ of genus~1. Moreover, $C$ has a
$k$-rational point $Q$ (the residual point to $P$ in $X\cap\Pi$),
so $(C,Q)$ is an elliptic curve.
The map $j_P$~acts by $-1$ in the group
law on this elliptic curve.

This clearly defines an involution of~$X$ and it can be shown
that $j_P$ is a rational, and hence birational, map.  Moreover,
if we suppose that $X$ is minimal, following Iskovskikh the
linear system corresponding to~$j_P$ is~$\left|(2d-1)A - 2dP\right|$.
As above, minimality is sufficient, but not necessary, for this linear system
to give the right map.

There is an extra detail when $d=2$: to obtain a map to 
$\P^3(1,1,1,2)$,
we must also compute sections of $\left|(4d-2)A -
  4dP\right|$ and find an additional section not in the
subspace spanned by products from $\left|(2d-1)A - 2dP\right|$.

The Bertini involution is described as a link of type~II in the
same way as the Geiser above; the only difference is that the
contracting curve $F \subset Y$, where $Y$ is the
blowup of $P \in X$, is now the unique effective curve in
$|(2d-2)A - (2d-1)E|$.

\bigskip Both the Geiser and Bertini involutions are naturally
defined algebraically using the pluri-anticanonical model of~$Y$.
This description clarifies many of the points we have touched
on---for example, it explains why it is natural that the Bertini
involution should act by $-1$ in the group law on the elliptic
curve $(C,Q)$ above---and it generalises to higher dimensions.
See~\cite{CPR} for explanation and many examples.

\section{Imposing conditions on functions}\label{sec:2}
This is the main algorithm. Let $X\subset\P^d$ be a surface and $P\in X$
a point ($k$-rational or not).
Let $\cH \subset |nA|$ be a linear system on $X$
that is the restriction of some linear system on $\P^d$.
We want to compute the linear system $\cH(-mP)\subset \cH$ of divisors
that have multiplicity at least $m\in\N$ at $P$.

There is one point to note. In solving these equations we compute a system
of homogeneous polynomials on the ambient space $\P^d$. Although we can,
and do, use these to make a map $\P^d\broken\P^d$, we only work with
the restriction of this to $X$---the extension to $\P^d$ is not
determined by the theory.

\subsection{Basic case: $\cH=|nA|$ and $P\in X$ is $k$-rational}
For this discussion, we assume that $X$ is a surface and $P\in X$
a rational point, the case we need for our application.  However,
the description extends to the blowup of any nonsingular point on
a variety of any dimension. It also applies to blowups of
higher-dimensional centres $\Gamma\subset X$, as long as one can
compute in generic coordinates along~$\Gamma$.

\paragraph{Normal coordinates along the blowup of $P$.}
The idea is to compute the blowup $Y$ of $P\in X$ as a sequence
of implicit functions on a blowup patch of the ambient~$\P^d$.
The functions defining $X$ are denoted $F_1,\ldots,F_r$ (that
is, these form a basis for the ideal of~$X$).

Let $c=d-2$ be the codimension of $X\subset\P^d$, so $c \le r$.
It is easy to describe an affine patch $\phi\colon
\A^d\rightarrow\P^d$ of the blowup of $P\in \P^d$.  Furthermore,
we may assume we have coordinates $u_1,\dots,u_d$ on $\A^d$ for
which the exceptional divisor $E\subset Y$ of the blowup of $P
\in X$ is the $u_{c+1}$-axis ($= u_{d-1}$-axis) and $u_d$ is not a critical
direction for the functions~$F_j$: that is, $(\partial
F_j/\partial u_d)_{j=1,\ldots,r}\not=0$ along $E$.  Roughly, we
work in coordinates $u_{d-1},u_d$ on $Y$ in a neighbourhood
of the generic point of the exceptional divisor~$E$; precisely,
this means working over the function field $k(E)=k(u_{d-1})$ and
regarding $u_d$ as a formal power series variable.  With these
variables as coefficients, the equations defining $Y$ are
polynomials in $u_1,\dots,u_c$, and there is a basis
$f_1,\dots,f_c$ of polynomials of this ideal.  We may assume, by
using pairwise resultants if necessary, that these polynomials
are each univariate in one of the variables: that is, $f_i =
f_i(u_i)$ for each $i=1,\dots,c$.  Now the implicit function
theorem guarantees the existence of power series $\phi_i(u_d)$
over~$k(E)$ for which $f_i(\phi_i)=0$.  These are computable to
any given precision; $u_i=\phi_i$ are the equations of~$Y$, and
we can use them to eliminate the variables $u_1,\dots,u_c$.

\paragraph{Compute $\cH$ along the exceptional divisor.}
Let $h$ be the generic element of~$\cH$, a homogeneous polynomial
of degree $n$ with unknown coefficients; we may assign
indeterminates $a_1,\dots,a_N$ as these coefficients, where $N$
is the number of monomials of degree $d$ on $\P^3$.  Pulling $h$
back to the normal coordinates to $E$ on $Y$ that we computed
above expresses it as a power series in $u_d$ with coefficients
in $k(E)$, say $q_0 + q_1u_d + q_2u_d^2 + \cdots$, computed to
precision at least $u_d^m$, with coefficients $q_i=q_i(u_{c+1})$
rational functions in $k(E)$ and in the unknown coefficients
$a_1,\dots,a_N$ of~$h$.  The latter only appear linearly in the
numerator of each $q_i$---pulling the expression $a_1x^d +
a_2x^{d-1}y+\cdots$ back to the normal coordinates evaluates the
monomials $x^d, x^{d-1}y,\dots$ at expressions in those
coordinates and then gathers terms together first in powers of
$u_d$ and then in powers of $u_{c+1}$, the variable along $E$.
The order of this power series is the degree of vanishing of
$\cH$ along $E$, so sections of $\cH(-mP)$ are those polynomials
whose coefficients are solutions of the first $m$ coefficients of
the power series, thought of as systems of linear equations in
$a_1,\dots,a_N$.

These systems can be solved using standard computational algebra
tools, and a basis for the solution space provides the
coefficients of homogeneous polynomials of degree $n$ that base
$\cH(-mP)$.

\subsection{Modifications for more general cases}
In general, we do not work with the full linear system $|nA|$ but
with some subsystem $\cH\subset |nA|$. This makes no essential
difference: we use a basis of the sections of~$\cH$ in place of
the monomials of degree $n$.

The main routine must also impose conditions at points $P\in X$
of higher degree.  We may assume that $P$ is irreducible over
$k$; it is enough then to make a finite field extension $k\subset
k_1$ that splits $P$ into $k_1$-geometric points and to apply the
algorithm at just one of these.  This determines $k_1$-linear
conditions on the coefficients of the linear system, which in
turn determine $\deg P$ linear conditions over~$k$.

\section{Del Pezzo surfaces of low degree}\label{sec:3}
To construct links of type~II on del Pezzo surfaces we must
compute Geiser and Bertini involutions given a centre $P\in X$.
We concentrate here on surfaces $X$ of degree $3$ or~$4$; we
discuss surfaces of degrees~2 and~5 in Section~\ref{sec:future}
below.

\subsection{Applying the algorithm}
Suppose $P\in X$ is an irreducible closed point.  The calculation
is in two steps.  If $P$ has degree $d-2$, compute a basis for
the sections of $|(d-1)A -dP|$ to make a Geiser involution; if
$P$ has degree $d-1$, compute a basis for the sections of
$|(2d-1)A -2dP|$ to make a Bertini involution.  In either case,
this basis has $d+1$ elements and defines a map $\psi \colon
\P^d\broken \P^d$, which (when restricted to~$X$) is the required
involution up to the choice of basis we made for the linear
system.  But without the correct basis, $\psi$ is unlikely even
to map $X$ to itself.  There is a unique $k$-linear automorphism
$\theta$ of~$\P^d$ such that $\theta\circ\psi$ is the required
involution; we must find this `missing' automorphism~$\theta$.

To determine $\theta$ we use the geometric definition of the
involution.  The aim is to find $d+2$ points of~$X$ that
span~$\P^d$ and whose images can be computed; these points,
together with their images, are enough to determine $\theta$. As
usual, the existence of many $k$-rational points is not expected,
but intersecting $X$ with random linear spaces containing $P$
provides arbitrary numbers of closed points of low degree.

Specifically, in the Geiser case we take general linear spaces of
dimension $d-2$ through $P$ and compute the residual intersection
with $X$ (using an ideal quotient to remove $P$ from the
intersection). Over the closure $\kbar$, this consists of two
points $q_1,q_2\in \Xbar$. We can compute $Q_1=\psi(q_1)$ and the
missing automorphism $\theta$ must satisfy $\theta(Q_1) = R_1$,
where we set $R_1=q_2$.  In the Bertini case we take hyperplanes
$\Pi$ through $P$ and intersect with~$X$. If $\Pi$ is
sufficiently general then $X\cap \Pi$ is nonsingular and we can
compute a Weierstrass normal form over~$k$; general fibres of the
elliptic involution are defined over~$k$, so they give points of
degree~2 on $X$ whose components $q_1,q_2\in \Xbar$ are exchanged
by the Bertini involution.  Again we can compute $Q_1 =
\psi(q_1)$ and the missing automorphism $\theta$ must satisfy
$\theta(Q_1) = R_1$, where we set $R_1=q_2$.

It remains to solve for $\theta$ given $d+2$ pairs of independent
points $Q_i,R_i$ as above.

\subsection{Finding the missing automorphism} 
Let $Q_1,\dots,Q_{d+2}$ and $R_1,\dots,R_{d+2}$ be two sequences
of points of~$\P_\kbar^d$, each of which spans~$\P_\kbar^d$, and
let $\theta$ be the $\kbar$-linear automorphism taking each $Q_i$
to $R_i$.  Standard linear algebra routines compute a matrix
representing~$\theta$, but in our application there are two
additional aspects to consider. First, our $\theta$ is defined
over~$k$ so the corresponding matrix should have entries in~$k$
too; by itself, this is not a problem.  Second, in practice we do
not work over~$\kbar$ but over a different extension $k\subset
k_i\subset\kbar$ for each pair of points $Q_i$ and~$R_i$; this
makes it difficult to apply the solution algorithm for linear
equations directly.  To deal with this we rephrase the algorithm
slightly and compute a matrix representing $\theta$ as follows.

Let $M_Q$ be a matrix whose rows are representatives for
$Q_1,\dots,Q_{d+2}$ and~$M_R$ a similar matrix for the $R_j$.  We
seek a $(d+1)\times (d+1)$ matrix $M$ for which $M_Q M = M_R$.
It is enough to solve this up to scalar multiples of rows; so,
letting $K_{R_i}$ be a matrix whose columns are a basis of $\ker
R_i$, it is enough to solve the system of equations
\begin{equation}\label{eq!ilin}
Q_1 M K_{R_1}=0, \quad\dots,\quad Q_{d+2}MK_{R_{d+2}} = 0.
\end{equation}
Each of these imposes $d$ linear conditions on $M$.  Since the
equations are independent, the solution space has dimension
$(d+1)^2 - d(d+2) = 1$, and the entries of~$M$ (itself only
defined up to a scalar) are then the coefficients of any
nontrivial solution.  This does not yet solve the problem: the
coefficients appearing in \eqref{eq!ilin} still lie in the
various fields~$k_i$ and it would be expensive to compute a
composite field containing them all.

Taking the trace of each of the equations gives a new system
\begin{equation}\label{eq!ilintrace}
\Tr_{k_1/k} (Q_1 M K_{R_1})=0, \quad\dots,\quad
\Tr_{k_{d+2}/k}(Q_{d+2}MK_{R_{d+2}}) = 0
\end{equation}
defined over~$k$. The following lemma is elementary; the main
point is to avoid $k$-linear relations holding between the chosen
representatives of the $Q_i$.

\begin{lem}
  Let $Q_1,\dots,Q_{d+2}$ and $R_1,\dots,R_{d+2}$ be two
  sequences of points of~$\P_\kbar^d$ as above. For fixed
  representatives of the $Q_i$ and~$R_i$, consider the systems of
  linear equations \eqref{eq!ilin} and~\eqref{eq!ilintrace}
  defined on $V=k^{(d+2)^2}$ as above.
  Denote the solution sets of these two systems
  by $W_{\eqref{eq!ilin}}$ and $W_{\eqref{eq!ilintrace}}$
  respectively.

  Then $W_{\eqref{eq!ilin}} \subset W_{\eqref{eq!ilintrace}}$ and
  equality holds for a general choice of representative for
  each~$Q_i$.
\end{lem}

\subsection{Degree 4 del Pezzo surfaces}
To illustrate our algorithms, we give examples to show how to
construct involutions using an implementation in \magma;
the implementation, together with these examples and their output,
is available at \cite{BKR}.

First make a nonsingular surface
\[
X\colon \left(\begin{array}{c}
xy - zt + 2x^2 + s^2 = 0\\
-x^2 + y^2 - z^2 + t^2 - s^2 = 0
\end{array}\right)\subset\P^4
\]
defined over~$k=\Q$.
{\small
\begin{verbatim}
> P4<x,y,z,t,s> := ProjectiveSpace(Rationals(),4);
> f := x*y - z*t + 2*x^2 + s^2;
> g := -x^2 + y^2 - z^2 + t^2 - s^2;
> X := Scheme(P4, [f,g]);
\end{verbatim}
}\noindent
For a Geiser involution we need a point $P \in X$ of degree~2,
which we construct by intersecting $X$ with a particular line:
{\small
\begin{verbatim}
> P := Intersection(X, Scheme(P4, [x,z,s]));
> Degree(P);
2
\end{verbatim}
}\noindent
This is all the data needed to construct the Geiser involution
centred in $P$. The map is stored in \magma\ as a composition of simpler
maps, so we use {\tt Expand(G)} to see its defining equations
of degree $d-1=3$; it is usually costly to do this step and
is unnecessary.
We also check that the involution really does map $X$ to itself.
{\small
\begin{verbatim}
> G := GeiserInvolution(X, P);
> Expand(G);
Mapping from: Prj: P4 to Prj: P4
with equations :
4/3*x*z^2 + 2/3*x*z*t - 1/3*y*z*t - 1/3*x*t^2 - 1/3*x*s^2
    + 1/3*y*s^2
-2/3*x*z^2 - y*z^2 - 7/3*x*z*t + 2/3*y*z*t + 2/3*x*t^2
    - 1/3*x*s^2 - 2/3*y*s^2
y^2*z + z*t^2 - z*s^2
4*y^2*z - 4*z^3 - y^2*t + 4*z*t^2 - t^3 - 2*z*s^2 + t*s^2
y^2*s + t^2*s - s^3
> G(X) eq X;
true
\end{verbatim}
}\noindent

For a Bertini involution we need a point $Q \in X$ of degree~3,
which we construct as the residual intersection to a $2$-plane
$\Pi$ containing a rational point $(0,1,1,0,0) \in X$.
The {\tt Support(Z)} command below computes the $k$-rational
support of~$Z$.
{\small
\begin{verbatim}
> Pi := Scheme(P4, [x+y-z,s]);
> Z := Intersection(X, Pi);
> supp := Support(Z); supp;
{ (0 : 1 : 1 : 0 : 0) }
> Degree(Z);
4
> R := Cluster(Representative(supp));
> Q := Scheme( P4, ColonIdeal(Ideal(Z), Ideal(R)) );
> B := BertiniInvolution(X,Q);
\end{verbatim}
}\noindent
Computing the map $B$ takes about half a minute, giving a
map defined by polynomials of degree $2d-1=7$; the equations of
the map are large, and so here we only show the initial terms.
{\small
\begin{verbatim}
> B;
Mapping from: Prj: P4 to Prj: P4
Composition of Mapping from: Prj: P4 to Prj: P4
with equations :
y^2*z^5 - 301231/288*y^2*z^4*t - 102767/144*x*z^5*t
    - 11755/32*y*z^5*t + 101059/72*z^6*t
    - 4791269/4608*y^2*z^3*t^2 - 2205985/2304*x*z^4*t^2 - ...
> B(X) eq X;
true
\end{verbatim}
}\noindent
These maps are indeed involutions and it is not necessary to
check this explicitly---although one can check, for instance,
that $G\circ G$ is the identity map by using the interpolation
routines explained below.

\subsection{Examples on a cubic surface}\label{sec:X3}
First make the surface $X\colon (x^3 + 2y^3 + 3z^3 + 4t^3=0)\subset\P^3$
defined over~$k=\Q$.
{\small
\begin{verbatim}
> P3<x,y,z,t> := ProjectiveSpace(Rationals(), 3);
> X := Scheme(P3, x^3 + 2*y^3 + 3*z^3 + 4*t^3);
\end{verbatim}
}\noindent
To make a Geiser involution, we need to choose a rational point
for its centre. We do not have code for finding rational points---in general,
this is an unsolved problem---but in this case we can see some obvious
choices: $(1,-1,-1,1)\in X$, for instance.
{\small
\begin{verbatim}
> p := X ! [1,-1,-1,1];
> G := GeiserInvolution(X,p);
\end{verbatim}
}\noindent
For a Bertini involution, we need a point of degree~2.
The line $L = (y+z+t = x-z+t=0)$ meets $X$ in the rational point
$(3,1,1,-2)$, and the residual intersection to this is an
irreducible point of degree~2.
{\small
\begin{verbatim}
> L := Scheme( X, [y+z+t, x-z+t] );
> Z := [ Y : Y in IrreducibleComponents(L) | Degree(Y) eq 2 ][1]; 
> B := BertiniInvolution(X,Z);
\end{verbatim}
}\noindent
(The Bertini calculation takes around 4 seconds; by comparison,
the Geiser calculation is instant.)
The map $B$ is defined by fairly large polynomials of degree~5.

We can construct another birational selfmap by composing these.
{\small
\begin{verbatim}
> h := Expand(B * G); h;
Mapping from: Prj: P3 to Prj: P3
with equations : 
54518131/19784704*x^4*y^6 + 59844679/7419264*x^3*y^7 + 
    272174051/29677056*x^2*y^8 + 7725367/1391112*x*y^9 +
    9154681/5564448*y^10 + 82450383/9892352*x^4*y^5*z + ...
\end{verbatim}
}\noindent
The equations of~$h$ are large polynomials of degree~10
that run over several pages.

\section{The group of birational selfmaps of a cubic surface}\label{sec:group}
It is well known, \cite{M} Theorem 38.1 for instance, that Geiser
and Bertini involutions (together with the subgroup of linear
automorphisms) generate the group $\Bir(X)$ of birational
selfmaps of any minimal cubic surface $X=X_3\subset\P^3$.  This
also follows from Iskovskikh's classification of links \cite{I}
Theorem~2.6: the only elementary links from cubic surfaces are
birational selfmaps.  (This is not the case in degrees 4 and 5
where factorisations of birational selfmaps into elementary links
may pass through other surfaces.)

The proof of Theorem~\ref{thm:decomp} works by induction on the
degree of the given selfmap~$\phi$.  The idea is to find a
basepoint of degree~1 or 2 that has high multiplicity in curves
belonging to the linear system defining~$\phi$.  Given such a
point, a so-called {\em maximal centre}, the proof precomposes
$\phi$ by the Geiser or Bertini involution it determines; this
decreases the degree and the induction continues.  We implement
this algorithmic step, demonstrating it here by factorising the
map $h$ computed in Section~\ref{sec:X3}.

First find the base locus of~$h$.
This is done naively, setting the defining equations of
$h$ to be zero. The result could be strictly bigger than
the base locus, but since we check the multiplicity of base
points later this does not matter.
{\small
\begin{verbatim}
> base_h := Scheme(X,DefiningEquations(h));
> Dimension(base_h);
0
> Degree(base_h);
205
\end{verbatim}
}\noindent
We need to identify irreducible components of this
base locus, and this could be a problem if there
really were 205 base points. But of course the
calculation above has found a non-reduced scheme,
and since we only need to know the base points set theoretically,
we can reduce it before further analysis.
Unfortunately, this seems to be difficult;
it takes about 4 minutes.
{\small
\begin{verbatim}
> base_red := ReducedSubscheme(base_h);
> Degree(base_red);
3
\end{verbatim}
}\noindent
Since the map $h$ is not linear, there must be a maximal centre
in this base locus. Untwisting by Bertini involutions is
likely to reduce the degree of~$h$ more dramatically than
by Geiser involutions, so we look for maximal centres
of degree~2 first.
{\small
\begin{verbatim}
> Q := [ Y : Y in IrreducibleComponents(base_red) | Degree(Y) eq 2 ][1];
\end{verbatim}
}\noindent
Of course, $Q$ is exactly the enforced base point $Z$ from
section~\ref{sec:X3} above, but without knowing that we need to check
that it is a maximal centre.
{\small
\begin{verbatim}
> assert is_maximal_centre(X, Q, h);
\end{verbatim}
}\noindent
This takes about a minute: it runs the main algorithm in high degree
to check that $Q$ has multiplicity strictly greater than 10 in
the linear system defining $h$.
We could omit this step since checking that the degree of~$h$ is reduced
after untwisting is sufficient.

We must now untwist $h$ by the Bertini involution centred in $Q$.
(Notice the order of composition in \magma: this is $h_1=h\circ\eps_1$.)
{\small
\begin{verbatim}
> eps1 := BertiniInvolution(X, Q);
> h1 := eps1 * h;
\end{verbatim}
}\noindent
There is a practical computational point here: if we expanded
out the equations of~$h_1$, they would have degree $10\times 5$,
the product of the degrees of the equations of~$\eps_1$ and~$h$.
But untwisting is meant to reduce the degree.
Working modulo the equation of~$X$, the equations of~$h_1$
have a large common factor that can be cancelled, but it is
not clear how to do this calculation.

Instead we use interpolation to find the the correct equations.
Bertini involutions reduce the degree by $4(\mult_Q(h) -
\deg(h))$, so even without knowing $\mult_Q(h)$, the possible
degrees for the resulting equations are limited.  Given a target
degree, the interpolation evaluates $h_1$ at many points of~$X$
and uses this collection of points and their images to impose
linear conditions on the coefficients of the desired equations.
As usual, we do not have a supply of rational points of~$X$ to
work with (there may be none), but we can intersect $X$ with
random rational lines to get points of degree~3 and use these.
If the target degree was too low, the solution space will only
include multiples of the defining equation of~$X$.  The first
time there are other solutions, these will be the coefficients of
the map.  The function we use below returns a boolean value,
which is false unless there is a unique additional solution.  If
the boolean value is true, the function also returns the
equations of~$h_1$. (The {\tt assert bool} statement causes a
crash unless the boolean is true.)  {\small
\begin{verbatim}
> bool,eqns := interpolate(X,h1,2);
> assert bool;
\end{verbatim}
}\noindent
We rebuild $h_1$ with the lower-degree equations.
{\small
\begin{verbatim}
> h1 := map< P3 -> P3 | eqns >; h1;
Mapping from: Prj: P3 to Prj: P3
with equations : 
x*y + y^2 + 3/2*x*z + 3/2*z^2 + 2*x*t - 2*t^2
1/2*x^2 + 1/2*x*y + 3/2*y*z - 3/2*z^2 + 2*y*t + 2*t^2
1/2*x^2 - y^2 + 1/2*x*z + y*z + 2*z*t + 2*t^2
-1/2*x^2 + y^2 + 3/2*z^2 + 1/2*x*t + y*t + 3/2*z*t
\end{verbatim}
}\noindent

Repeat the process with $h_1$:
{\small
\begin{verbatim}
> base_h1 := Scheme(X, DefiningEquations(h1));
> Dimension(base_h1);
0
> base1_red := ReducedSubscheme(base_h1);
> Degree(base1_red);
1
> r := Representative(Support(base1_red)); r;
(1 : -1 : -1 : 1)
\end{verbatim}
}\noindent
So the base locus is a single rational point $(1, -1, -1, 1)\in X$.
There is no choice but to untwist by the Geiser involution.
{\small
\begin{verbatim}
> eps2 := GeiserInvolution(X, r);
> h2 := eps2 * h1;
> bool,eqns := interpolate(X,h2,1);
> assert bool;
> h2 := map< P3 -> P3 | eqns >; h2;
Mapping from: Prj: P3 to Prj: P3
with equations : x, y, z, t
\end{verbatim}
}\noindent
The resulting map $h_2$ gives a linear automorphism of~$X$
so the factorisation is complete.
In this case, it is easy to check that the identity is the only
automorphism of~$X$, so the group $\Bir(X)$ is generated
by all Geiser and Bertini involutions---what this group is,
therefore, has become an arithmetic question involving the low degree
points of~$X$ and the relations between Geiser and Bertini involutions.

\section{Rational surfaces and \magma}\label{sec:future}
We describe briefly the broader context of a possible complete
implementation of minimal rational surfaces in \magma\ that builds on
the algorithms here.  It is realistic to expect: to compute the
factorisation of rational maps between two-dimensional Mori fibre
spaces, the so-called {\em Sarkisov programme};
to address rationality questions; and to analyse elliptic fibrations.
Beck and Schicho have algorithms, implemented in \magma,
that compute the Picard group of a del Pezzo surface and, if necessary,
carry out the minimal model programme in certain circumstances.

\subsection{The Sarkisov programme for rational surfaces}
We were able to carry out the Sarkisov programme for selfmaps of
a cubic surface in Section~\ref{sec:group} because in that case
other surfaces are not involved.  The same is true for del Pezzo
surfaces of degrees~1 and 2, but in other cases we need the full
classification of links in order to proceed.  An implementation
of the full Sarkisov programme for factorising birational maps
between minimal rational surfaces, following Iskovskikh---so
realising Theorem~\ref{thm:decomp} explicitly---is feasible,
given the following additional components.

The involutions on del Pezzo surfaces of degrees~2 and 5 are
essentially the same as those described here. In degree~5 the
calculations are currently too slow to be practical.  In degree~2
we must work with weighted projective space; this only
complicates the calculation slightly. The main additional
difficulty is in recognising the surface as an image of the
multi-linear system that determines the involution. For example,
if $P=(1,0,0,0)\in X\subset\P(1,1,1,2)$ with tangent plane
$T_P(X) = (y=0)$, in coordinates $x,y,z,t$, then
\[
|3A - 4P| = \left< y^3, y^2z, yz^2, y^2x \right>
\]
and $|6A - 8P|$ is spanned by quadratic expressions in these
together with $y^4t$ and~$y^3zt$.
These determine a map $X\broken \P^5(1^4,2^2)$ whose
image $Y$ is isomorphic to $X$---but we would need to compute
$|-K_Y|$ and $|-2K_Y|$ to make that identification (and then
study the geometrical description of the Bertini involution as
before to make the right choice of automorphism).

Computing automorphisms of del Pezzo surfaces is straightforward.
Since del Pezzo surfaces are embedded by the anticanonical class,
any automorphism $X\rightarrow X$ extends to an automorphism of
the ambient projective space. If this space is $\P^d$ then the
only automorphisms are projective linear maps; when
$X=X_4\subset\wP^3=\P^3(1,1,1,2)$ we must allow quasi-linear
maps as follows:
\[ 
x_1\mapsto f_1(x_1,x_2,x_3),
\dots,\quad
x_3\mapsto f_3(x_1,x_2,x_3),
\quad
y\mapsto ay + g(x_1,x_2,x_3)
\]
where $x_1,x_2,x_3,y$ are the homogeneous coordinates on~$\wP^3$,
the $f_i$ are linear forms, $g$~is a quadratic form and $a\in k$
does not have a square root in~$k$.

General conic bundles can be embedded in scrolls, even if the base
rational curve does not have a rational point. With such a description,
one can compute links by explicit equations as we
have for del Pezzo surfaces.
Cremona and van Hoeij \cite{CvH} give an algorithm for
Tsen's theorem (over an extension of~$k$, as necessary),
which is one of the essential tools for this.

\subsection{Rationality questions for minimal rational surfaces}\label{sec:rat}
While we assume that our surfaces $X$ are rational over~$\kbar$,
the question of whether they are also rational over~$k$ is well
studied.  Swinnerton-Dyer \cite{SD} answers the question for a
cubic surface in terms of a rational point on $X$ and conditions
on the configuration of the 27 lines.  In degree~5, Enriques
proved that $X$ automatically has a rational point and is
rational over~$k$; see~\cite{Sh-B}.  In degrees~$\ge 6$, $X$ is
rational over~$k$ if and only if it has a rational point, and
either case can happen.  In degrees 8 and 9, de
Graaf--Harrison--Pilnikova--Schicho \cite{dG} determine whether a
del Pezzo surface (not necessarily minimal, and for $k$ a number
field) is rational and they compute a parametrisation over~$k$ in
that case; this is implemented in \magma.

\subsection{Elliptic fibrations on minimal rational surfaces}
Elliptic fibrations birational to minimal rational surfaces are
classified in some cases: Dolgachev \cite{D} classifies elliptic
fibrations birational to~$\P^2$ (the construction of these
fibrations dates back to Halphen~\cite{H}); Cheltsov~\cite{Ch}
and Brown--Ryder~\cite{BR} analyse elliptic fibrations birational
to minimal cubic surfaces; and Cheltsov also analyses the case of
del Pezzo surfaces of degrees~1 and~2.  A birational map from $X$
to an elliptic fibration can be regarded as a limiting case of a
birational map to another Mori fibre space, so the methods of
construction and exclusion for the two problems are very similar.
The computational aspects of the elliptic fibration problem are
fully analysed in the degree~3 case in~\cite{BR}, and a full
\magma\ implementation is given; this is not yet done in the
other cases.


\vspace{5mm}
\noindent
Gavin Brown, IMSAS, University of Kent, Canterbury, CT2 7AF, UK.\\
{\small \texttt{gdb@kent.ac.uk} }

\vspace{5mm}
\noindent
Alexander Kasprzyk, IMSAS, University of Kent, Canterbury, CT2 7AF, UK.\\
{\small \texttt{A.M.Kasprzyk@kent.ac.uk} }

\vspace{5mm}
\noindent
Daniel Ryder, Department of Mathematics, University Walk, Bristol, BS8 1TW, UK.\\
{\small \texttt{daniel.ryder@bristol.ac.uk} }
\end{document}